\documentclass{dcds}
\usepackage{amsfonts}

\newcommand{\N}{\Bbb N}
\newcommand{\R}{\Bbb R}

\def\qed{\hbox to 0pt{}\hfill$\rlap{$\sqcap$}\sqcup$}
\newtheorem{theorem}{Theorem}[section]
\newtheorem{proposition}[theorem]{Proposition}
 \newtheorem{lemma}[theorem]{Lemma}

\theoremstyle{definition}

\theoremstyle{remark}
\newtheorem{remark}[theorem]{Remark}

\numberwithin{equation}{section}
\begin{document}


\title[Admissible wavefront speeds]{Admissible wavefront speeds for a single species
reaction-diffusion equation with delay}

\author[Elena Trofimchuk and  Sergei Trofimchuk]{}
\email{trofimch@imath.kiev.ua}
\email{trofimch@inst-mat.utalca.cl}


\date{\today}

\maketitle

\centerline{\scshape  Elena Trofimchuk}
\medskip

{\footnotesize \centerline{ Department of Mathematics II, }
\centerline{ National Technical University `KPI', Kyiv,  Ukraine}
}
\medskip

\centerline{\scshape  Sergei Trofimchuk}
\medskip

{\footnotesize \centerline{Instituto de Matem\'atica y
F\'{\i}sica} \centerline{ Universidad de Talca, Casilla 747,
Talca, Chile } }

\bigskip
\begin{quote}{\normalfont\fontsize{8}{10}\selectfont
{\bfseries Abstract.} We consider equation $u_t(t,x) = \Delta
u(t,x)- u(t,x) + g(u(t-h,x))\ (*) $, when  $g:\R_+\to \R_+$ has
exactly two fixed points: $x_1= 0$ and $x_2=\kappa>0$. Assuming
that $g$ is unimodal and has negative Schwarzian, we indicate
explicitly a closed interval $\mathcal{C} =
\mathcal{C}(h,g'(0),g'(\kappa)) = [c_*,c^*]$ such that $(*)$ has
at least one (possibly, nonmonotone) travelling front propagating
at velocity $c$ for every $c \in \mathcal{C}$. Here $c_*>0$ is
finite and $c^* \in \R_+ \cup \{+\infty\}$. Every time when
$\mathcal{C}$ is not empty, the minimal bound $c_*$ is sharp so
that there are not wavefronts moving with speed $c < c_*$. In
contrast to reported results, the interval $\mathcal{C}$ can be
compact, and we conjecture that some of equations $(*)$ can indeed
have an upper bound for propagation speeds of travelling fronts.
As particular cases, Eq. $(*)$ includes the diffusive Nicholson's
blowflies equation and the Mackey-Glass equation with nonmonotone
nonlinearity.
 \par}
\end{quote}


\section{Introduction}
\label{sec.int} \noindent \label{intro} In this paper, we study
the existence of positive and generally nonmonotone travelling
waves for a family of delayed reaction-diffusion equations
\begin{equation}\label{17}
u_t(t,x) = \Delta u(t,x)  - u(t,x) + g(u(t-h,x)), \ u(t,x) \geq
0,\ x \in \R^m,
\end{equation}
which  has exactly two non-negative equilibria $u_1 \equiv 0, \
u_2 \equiv \kappa >0$ (so that $g(\kappa) = \kappa, \ g(0) =0$).
The nonlinearity $g$ is called {\it the birth function}. It is
therefore non-negative, and generally nonmonotone and bounded. In
particular, with $g(u)= pue^{-u}$ in (\ref{17}), we get the
diffusive Nicholson's blowflies equation proposed in \cite{sy}.
Since the biological interpretation of $u$ is the size of an adult
population, we will consider {\it only} nonnegative solutions for
(\ref{17}).

Recently, the existence of  travelling fronts connecting the
trivial and positive steady states in (\ref{17}) (and its
non-local multi-dimensional generalizations) was studied in
\cite{fhw,FT,ma,ma1,swz,wz} by means of two essentially different
approaches (see also \cite{fhw,gouk,gouS, gouss} for  other
methods which can be applied to analyze this problem). The first
method, which was proposed in \cite{sz,wz,zw} and then further
substantially developed in \cite{hz,hz1,ma,ma1,swz}, uses
positivity and monotonicity properties of the integral operator
\begin{equation*}\label{psea} \hspace{-7mm}
(Ax)(t) = \frac{1}{\epsilon (\mu -
\lambda)}\left\{\int_{-\infty}^te^{\lambda (t-s)}g(x(s-h))ds +
\int_t^{+\infty}e^{\mu (t-s)}g(x(s-h))ds \right\},
\end{equation*}
where $\lambda < 0 < \mu$ satisfy $\epsilon z^2 -z -1 =0$ and
$\epsilon^{-1/2} = c > 0$ is the front velocity. As it can be easily
observed, the profiles $x \in C(\R,\R_+)$ of travelling waves are
completely determined by the integral equation $Ax=x$ and the first
approach consists in the use of an appropriate fixed point theorem
to $A:K \to K$, where convex and closed set $K \subseteq \{x: 0\leq
\phi^-(t) \leq x(t) \leq \phi^+(t)\} \subset C(\R,\R_+)$ belongs to
an adequate Banach space $(C(\R,\R_+), |\cdot|)$. Hence, the
existence of travelling fronts can be established if we will be able
to realize  {\it Ma-Wu-Zou reduction} which consists in choosing
profile's subset $K$ of an appropriate Banach space $(C(\R,\R_+),
|\cdot|)$. This space should be 'nice' enough to assure the
compactness (or monotonicity) of $A$ while $\phi^{\pm}$ have to
satisfy $A\phi^- \geq \phi^-$ and $A\phi^+ \leq \phi^+$. All these
requirements are not easy to satisfy, and this why only relatively
narrow subclasses of (\ref{17}) (e.g. monotone or quasi-monotone in
the sense of \cite{wz}) were considered within this approach. Here,
we would like to mention the recent works \cite{hz,hz1,ma1,swz}
which open new interesting possibilities of this method. It should
be observed that a correct choice of $\phi^{\pm}$ is related to the
calculation of the minimal speed of wavefront's propagation, e.g.
see \cite{gouss,ma1,swz,sz}

The second method was proposed in \cite{fhw}. It essentially
relies on the fact that, in a 'good' Banach space, the Frechet
derivative of $\lim_{\epsilon \to 0} A$ along a heteroclinic
solution $\psi$ of the limit delay differential equation
$x'(t)=-x(t)+ g(x(t-h))$ is a surjective Fredholm operator. In
consequence, the Lyapunov-Schmidt reduction can be used to prove
the existence of a smooth family of travelling fronts in some
neighborhood of $\psi$. As it was shown in \cite{FT} for the case
of (\ref{17}), this family contains positive solutions as well.
Furthermore, the method of \cite{fhw, FT} allows to make a
conclusion about the uniqueness (up to translations) of the
positive solution to Eq. (\ref{17}) for every sufficiently small
$\epsilon$. Additionally, the following condition established in
\cite{FT}
\begin{equation} \label{gsc0}
e^{- h} > - \Gamma \ln \frac{{\Gamma}^{2} - \Gamma}{{\Gamma}^{2} +
1}, \quad \Gamma \stackrel{def}= g'(\kappa),
\end{equation}
and which is sufficient for the existence of positive travelling
fronts in (\ref{17}) for $g'(0) > 1$ and small $\epsilon
>0$, is also 'very close' to the sharp one (see
\cite{FT, ltt} for more details and references). However, this
analytical approach has one important limitation: it can be used
only to prove the existence of travelling fronts for sufficiently
small $\epsilon
>0$, or, in other words, for sufficiently large speeds $c >0$.
On the other hand, calculation of the minimal wave speed (and the
asymptotic speed of propagation) is very important in
applications, see discussion about the marginal stability
conjecture and related questions in \cite{BD,wim,tz,whw,Xin}.

This work is inspired by the above presented approaches and,
especially, by \cite{fhw,FT,FT1,gyt,ma,sz,wz}. For a broad family
of nonlinearities $g$,  we prove that Eq. (\ref{17}) has positive
travelling wave fronts $u(t,x) = \phi(ct + \nu\cdot x)$ provided
that $g'(0) >1$, and $c$ exceeds the minimal velocity: $c
> c_*(h,g'(0))$, and
\begin{equation} \label{gsc1}
\xi(h,c) = \frac{\mu - \lambda}{\mu e^{-\lambda h}- \lambda
e^{-\mu h}}
 \geq \frac{{\Gamma}^{2} + \Gamma}{{\Gamma}^{2} + 1}.
\end{equation}
These fronts generally are not monotone: in fact, they can
oscillate infinitely about the positive steady state. We also show
that, for large negative values of $s$, the wave profile $\phi(s,
c)$ is asymptotically equivalent to an increasing exponential
function. The minimal speed $c_*(h,g'(0))$ is determined as the
unique positive number for which the characteristic equation
$(z/c)^2-z-1+g'(0)\exp(-z h)=0$ has a real multiple root. This
value of $c_*(h,g'(0))$ is sharp: there are no travelling waves
moving at speed $c < c_*(h,g'(0))$.

Now, notice that $ 0 < e^{-h} \leq \xi(h,c)  < 1 $ and
$\xi(h,+\infty) = e^{-h}$. Thus inequality (\ref{gsc1}) can be
replaced by a stronger one:
\begin{equation} \label{gsc}
e^{- h} > \frac{{\Gamma}^{2} + \Gamma}{{\Gamma}^{2} + 1}.
\end{equation}
We show that, under the conditions $g'(0) > 1$ and (\ref{gsc}),
the speed interval $\mathcal{C}$ is infinite, in fact,
$\mathcal{C} = [c_*, + \infty)$. The latter result is an important
supplement to the existence conditions proved in \cite{FT} since
(\ref{gsc}) can viewed as a linearized version of (\ref{gsc0})
where $\ln(1+w)$ is replaced with $w \geq \ln(1+w), \ w \geq 0$.

 Next, if (\ref{gsc1}) is fulfilled and (\ref{gsc}) does
not hold, we still can prove the existence of a compact interval
$\mathcal{C} = [c_*, c^*]$ of admissible speeds. Here $c^*=
c^*(g'(\kappa),h)$ is determined from the equation
\begin{equation}\label{ce} \frac{\mu(c^*) - \lambda(c^*)}{\mu(c^*)
e^{-\lambda(c^*) h}- \lambda(c^*) e^{-\mu(c^*) h}}
=\frac{{\Gamma}^{2} + \Gamma}{{\Gamma}^{2} + 1}.
\end{equation}
In order to state the main result of this work, we need to present
our basic hypothesis:
\begin{description}
\item[{\rm \bf(H)}] Let $g \in C(\R_+, \R_+)$ have only one point
of local extremum $x_M$ (maximum) and assume that  $g(0)=0$, $g(x)
> 0$ if $x > 0$, and
$$1<
g'_-(0+)=\liminf_{x\to +0}g(x)/x \leq \limsup_{x\to +0}g(x)/x =
g'_+(0+) < +\infty.$$  We suppose further that the equation $g(x)=x$
has exactly two fixed points $0$ and $\kappa>0$. Set $\zeta_2 =
g(x_M)$ and let $\zeta_1 \leq \min \{g ^2(x_M),\ x_M, \kappa\}$ be
such that $g(\zeta_1) = \min_{s \in [\zeta_1,\zeta_2]}g(s)$.
Additionally,  $g \in C^3[\zeta_1,\zeta_2]$ and the Schwarz
derivative $Sg$ is negative for all $x \in
[\zeta_1,\zeta_2]$, $x \not=x_M$:\\
$(Sg)(x)=g'''(x)(g'(x))^{-1}-(3/2)
\left(g''(x)(g'(x))^{-1}\right)^2 < 0$.
\end{description}
In the sequel $c_*^+$ will denote the unique positive number for
which the characteristic equation $(z/c)^2-z-1+g'_+(0+)\exp(-z
h)=0$ has a real multiple root. Also, we set $\mathcal{C}_+ =
[c_*^+, c^*]$, where $c^*= \infty$ if $(\ref{gsc})$ holds, and
$c^*$ is determined from (\ref{ce}) if $(\ref{gsc})$ does not hold
and $(\ref{gsc1})$ is true for some positive $c$.
 Now we are ready to state our main
result:
\begin{theorem} \label{main} Assume {\rm \bf(H)}. Then, for every $c \in \mathcal{C_+}= [c_*^+, c^*]$ and
$\nu \in \R^m, \ \|\nu\| =1$, equation (\ref{17}) has a positive
wavefront $u(t,x) = \phi(ct + \nu\cdot x, c)$. Moreover,
$$c_*^+ \leq \min \left\{ 2\sqrt{g'_+(0+)-1}, \sqrt{\frac{\ln g'_+(0+)}{h}}\right\}$$
and  $c_*^+ = c_*^+(g'_+(0+),h)$  decreases monotonically to $0$ as
$ h \to +\infty$ for every fixed $g'_+(0+) > 1$. If additionally $g$
is $C^2$-smooth in some neighborhood of $x =0$, then the value of
$c_* = c_*^+$ is sharp (minimal speed of travelling waves). Finally,
the travelling profile $\phi(t)$ oscillates about $\kappa$ on every
interval $[\tau, +\infty)$ provided that the equation
\begin{equation}\label{cchh1} (z/c)^2-z-1+g'(\kappa)\exp(-z h)=0
\end{equation} does not have any root in $(-\infty,0)\cup i\R$.
\end{theorem}
We notice again that, in Theorem \ref{main}, the interval
$\mathcal{C}$ can be finite exhibiting a situation which is not
possible for the reaction-diffusion monostable equations without
delay \cite[Theorem 8.3 (ii)]{GK}  In fact, we conjecture that, for
the smooth nonlinearity $g$ satisfying {\rm \bf(H)}, it holds
$\mathcal{C} = [c_*, c_{opt}^*]$ where $c_{opt}^*=
c_{opt}^*(g'(\kappa),h)\in \R \cup \{\infty\}$ is defined as the
largest value for which characteristic equation (\ref{cchh1}) has
exactly one root in the open half plane  $\{\Re z > 0\}$. In
Proposition \ref{co}, we give an argument in favor of this
conjecture, establishing that $[c_*, c^*]\subseteq [c_*,
c_{opt}^*]$. The validity of the upper bound $c_{opt}^*$ will be
discussed in a forthcoming paper. Moreover, in view of the Smith
global stability conjecture \cite{ltt,smith} and Lemma \ref{c1} with
Remark \ref{rem} below, we believe that $\mathcal{C} = [c_*,
+\infty)$ if and only if  the scalar delay differential equation
\begin{equation}\label{17a}
u'(t) =   - u(t) + g(u(t-h)), \ u \geq 0,
\end{equation}
is globally stable.

The use of the Schwarz derivative in {\rm \bf(H)} is motivated by
the following proposition due to Singer, see \cite{gyt, LMT}:
\begin{proposition}\label{SW}
 Assume that $\sigma: [\zeta_*, \zeta^*] \to [\zeta_*, \zeta^*],
 \  \sigma \in C^3[a, b]$, is either strictly decreasing function or it has only one
critical point $x_M$ (maximum) in $[\zeta_*, \zeta^*]$. If the
unique fixed point $\kappa \in [\zeta_*, \zeta^*]$ of $\sigma$ is
locally asymptotically stable and the Schwarzian derivative
satisfies $(S \sigma)(x) < 0$ for all $x \not=x_M$ then $\kappa$
is globally asymptotically stable.
\end{proposition}
It is easy to see that Theorem \ref{main} applies to both the
Nicholson's blowflies equation and the Mackey-Glass equation with
nonmonotone nonlinearity \cite{ltt1}. On the other hand, the
assumption about the negativity of $Sg$ (which requires $C^3-$
smoothness of $g$) can be considerably relaxed with the use of a
generalized Yorke condition introduced in \cite{flot,ltt1,ltt}.

The structure of this paper is as follows: in the next section, we
study the dependence of roots to the characteristic equation
(\ref{cchh1}) in the positive steady state of (\ref{17}) on the
parameter $\epsilon = c^{-2}$. The same for characteristic equation
\begin{equation}\label{char3} (z/c)^2-z-1+g'(0)\exp(-z h)=0
\end{equation} in the trivial equilibrium. Remark \ref{rere} of this section
contains an estimation of the minimal speed $c_*$ and shows the
monotonicity of $c_*$ with respect to $h$. In Section 3, we
realize an appropriate Ma-Wu-Zou reduction, indicating a 'good'
Banach space $(C(\R,\R_+), |\cdot|)$ and its subset $K$. The
novelty of approach presented there consists in the very simple
form of the both bounds $\phi^{\pm}$ for $K$, we believe that our idea
can simplify selection of upper and lower solutions in other
similar situations. Section 3 contains our first results about the
existence of positive bounded solutions to the integral equation
$Ax =x$. These solutions are prototypes of wavefronts (in fact,
following the terminology of \cite{GK}, they are semi-wavefronts
to $0$). Before proving in Section 5, Theorem \ref{mainex}, that
they are indeed asymptotically constant at $+\infty$, we firstly
establish their persistence in Section 4. Furthermore, in Remark
\ref{rr} and Theorem \ref{ms} we prove that $c_*$ is indeed the
minimal travelling speed. Additionally, in Section 5 we are
concerned with the asymptotic formulae for travelling profile
$\phi(t)$ as $t \to -\infty$ (Remark \ref{rrra}) and discuss the
problem  of non-monotonicity of $\phi(t)$ as $t \to +\infty$
(Remark \ref{rrr}). Finally, in Proposition \ref{co} we show that
our main result is consistent with the above stated conjecture
about the exact speed interval for Eq. (\ref{17}).
\section{Characteristic equations}
In this section, we are concerned with the estimation of the
minimal speed of propagation of wavefronts.  For our model, this
speed is completely determined by the linearization of (\ref{17})
about the trivial steady state. The other question of interest is
the behavior of roots to the characteristic equations
(\ref{char3}) and (\ref{cchh1})
 as functions of parameters $\epsilon,
h$. Besides of the use of these results in Theorem \ref{main}, we
need them to give a partial justification of our conjectures
stated in Introduction. Everywhere below, we set
\begin{equation*}\label{char2}
\phi(z) =-z-1+a\exp(-z h), \quad \psi(z,\epsilon) =\epsilon
z^2-z-1+a\exp(-z h).
\end{equation*}
\begin{lemma} \label{c1} Take $a < 0$ and $h >0$ and  assume
that $\Re \lambda_j < 0$ for every solution $\lambda_j$ of the
characteristic equation $\phi(\lambda) =0$. Then, for every fixed
$\epsilon_0 > 0$ the equation $\psi(\lambda,\epsilon_0) =0$ has
only one root $\lambda_0$ in the right semi-plane $\Re z \geq 0$.
Moreover, $\lambda_0 \in (0,+\infty)$. \end{lemma}
\begin{proof}
The existence of the unique positive root $\lambda_0$ to
$\psi(\lambda,\epsilon_0) =0$ (with some fixed $\epsilon_0 >0$) is
evident. Next, every multiple root $z_0$ to this equation has to
satisfy the system
$$
\epsilon_0 z_0^2-z_0-1+a\exp(-z_0 h)=0, \ 2\epsilon_0
z_0-1-ah\exp(-z_0 h)=0,
$$
which implies
\begin{equation}\label{ep1}
(\epsilon_0 z_0^2-z_0-1)h+2\epsilon_0 z_0-1=0, \ a\exp(-z_0 h)=
\frac{2+z_0}{2+ hz_0}.
\end{equation}
The first equation of (\ref{ep1}) implies that $z_0$ is real while
the second equation of (\ref{ep1}) says that $z_0 < 0$. Now, it is
easy to see that $z_0$ is at most of the multiplicity 2, and that
  $(z_0, \epsilon_0)$ is a bifurcation
point where two real roots merge and disappear as $\epsilon \to
\epsilon_0 +$.

Hence, each root $\lambda_j \not\in \R$ of $\psi(z,\epsilon_0) =0$
determines a unique smooth function $\lambda_j(\cdot): (\alpha_j,
\beta_j) \to \mathbb{C}$ defined on some maximal open interval
$(\alpha_j, \beta_j)\subseteq [0, +\infty)$ containing $\epsilon_0$
and such that $\lambda_j(\epsilon_0) =\lambda_j,
\psi(\lambda_j(\epsilon), \epsilon) =0$. Additionally, we claim that
the arc $\lambda_j(\cdot):(\alpha_j, \beta_j) \to \mathbb{C} $ can
not intersect the imaginary axis from left to right. Indeed, we have
that
$$
\lambda_j'(\epsilon)= - \frac{z^2}{2\epsilon z -1 + h(\epsilon
z^2-z-1)},
$$
so that, at the moment $\epsilon_*$ of the eventual intersection we
have $\lambda_j(\epsilon_*) = i \omega$ and
\begin{equation}\label{der}
\Re \lambda_j'(\epsilon_*) = - \omega^2(1 +h +\epsilon
h\omega^2)/((1+h(\epsilon \omega^2 +1))^2+\omega^2(2\epsilon-h)^2)
<0,
\end{equation}
a contradiction.

Now, if $\lambda_j \not\in \R$ and $\lambda_j(\epsilon)$ were
bounded as $\epsilon \to \alpha_j+> 0$, then $\lambda_j(\alpha_j+) =
z_0$ would be the negative real root of the multiplicity 2.
In consequence, $\lambda_j(\epsilon)$ takes only real negative
values for $\epsilon>\alpha_j$, a contradiction. Therefore, if
$\alpha_j> 0$ and $\lambda_j \not\in \R$, then $\lim_{\epsilon \to
\alpha_j+}\Re(\lambda_j(\epsilon)) = - \infty$.

Next, in the event that $\alpha_j =0$ and
$\{\lambda_j(\epsilon_n)\}$ is bounded for some $\epsilon_n \to 0+$,
we can assume that $\lim \lambda_j(\epsilon_n)$ exists and satisfies
$\phi(\lambda) =0$. Therefore $\Re \lambda_j(\epsilon_n) < 0$
implying $\Re \lambda_j(\epsilon) < 0$ for all $\epsilon \in
(\alpha_j, \beta_j)$. If $\alpha_j =0$, and $\lambda_j(\epsilon)$ is
unbounded as $\epsilon \to 0+$, and $\Re(\lambda_j(\epsilon))
> 0$ in some right neighborhood of $0$, then
$\epsilon\lambda_j(\epsilon) \to 1$ as $\epsilon \to 0+$ so that
$\Re \lambda_j(\epsilon) \to +\infty$ when $\epsilon \to 0+$. It
should be noticed now that, for small $\epsilon
>0$ and sufficiently large $C >0$, the function $\psi(\lambda,\epsilon_0)$ has, in
virtue of the Rouch\'e theorem,  only one positive root in the half
plane $\Re z > C$. Hence, we can conclude that the unique branch
$\lambda_j(\epsilon)$ with infinite $\lambda_j(0+)$ corresponds to
the positive real root of $\psi(\lambda,\epsilon_0) =0$.
\end{proof}
\begin{remark} \label{rem} Analyzing proof of Lemma \ref{c1}, we see that if, for fixed
$a <0, h >0$ and for all $\epsilon > 0$,  equation
$\psi(\lambda,\epsilon_0) =0$ has a unique root in the half plane
$\{\Re z > 0\}$, then the quasipolinomial $\phi(\lambda)$ must be
stable. Furthermore, (\ref{der}) implies that, if for fixed $a <0,
h, \epsilon_0 > 0$, the equation $\psi(\lambda,\epsilon_0) =0$ has
a unique root in the half plane $\{\Re z
> 0\}$, then this property will be maintained  for all $\epsilon >
\epsilon_0$.
\end{remark}
Next proposition, which can be considered as an analog of Lemma
\ref{c1} for positive $a > 1$, is crucial for the calculation of
the minimal speed of propagation.
\begin{lemma} \label{L23} Assume that $a > 1$. Then there exists a positive $\epsilon_0 =
\epsilon_0(a,h)$ such that, for every fixed $\epsilon \in (0,
\epsilon_0)$, equation $\psi(\lambda,\epsilon) =0$ has exactly two
real roots $0 < \lambda_1(\epsilon) < \lambda_2(\epsilon)$. If
$\epsilon > \epsilon_0$, then $\psi(\lambda,\epsilon)> 0$ for all
$\lambda \in \R$. Moreover, the cardinal number $N(\epsilon)$ of
the set $\{\lambda_j: \psi(\lambda_j,\epsilon)=0 \ {\rm and } \
\Re \lambda_j
>0\}_{j=1}^{N(\epsilon)}$  is a decreasing function on $(0,\epsilon_0)$. Finally,
all roots of $\psi(\lambda,\epsilon) =0, \ \epsilon \in (0,
\epsilon_0)$ are simple and we can enumerate them in such a way
that
$$\dots \leq \Re
\lambda_4 = \Re \lambda_3 < \lambda_1 < \lambda_2.
$$
\end{lemma}
\begin{proof} The proof of the existence of such $\epsilon_0 = \epsilon_0(a,h)$ is
elementary and thus it is omitted here. In fact, $\epsilon_0$ can
determined from system (\ref{ep1}), where the second equation has
a unique positive solution $z_0 =z_0(a,h)$ for every $a >1$ and $h
>0$. Next, for every fixed $\epsilon >0$, the finiteness of the
set $\{\lambda_j: \Re \lambda_j >0\}_{j=1}^{N(\epsilon)}$ of all
roots with non-negative real parts to equation
$\psi(\lambda,\epsilon) =0$ is a well known fact (e.g. see
\cite[P.18]{hale}). The monotonicity of $N(\epsilon)$ can be
deduced from (\ref{der}) by repeating the same argument as in the
proof of Lemma \ref{c1}. The proof of simplicity of roots when
$\epsilon \in (0, \epsilon_0)$ is straightforward (e.g. see Lemma
13 in \cite{FT}). Finally, if $\epsilon > 0$ is sufficiently small
then the last assertion of Lemma \ref{L23} follows from Lemmas
7,13 of \cite{FT}. Take now the curves $\lambda_j(\epsilon)$ and
suppose for a moment that $\epsilon_1 \in (0, \epsilon_0)$ is the
minimal value for which $\Re \lambda_4(\epsilon_1) = \Re
\lambda_3(\epsilon_1) = \lambda_1(\epsilon_1) <
\lambda_2(\epsilon_1).$ Then the vertical line $\Re z =
\lambda_1(\epsilon_1)$ contains three simple roots of $\epsilon_1
z^2-z-1+a\exp(-z h)=0$. On the other hand, $\Re(\epsilon_1 z^2-z-1
+a\exp(-z h))=0$ with $\Re z = \lambda_1$ implies that
$$
ae^{-\lambda_1 h} \cos(h\Im z) = \epsilon(\Im z)^2 + (\lambda_1 +1
- \epsilon \lambda_1^2) = \epsilon(\Im z)^2 + ae^{-\lambda_1 h}.
$$
Therefore, we get
$$
2\epsilon(\Im z)^2 + ae^{-\lambda_1 h}\sin^2(0.5h\Im z) =0.
$$
Thus $\Im z =0$ and, in consequence, $\lambda_3(\epsilon_1) = \Re
\lambda_3(\epsilon_1) = \lambda_1(\epsilon_1)$, a contradiction.
\end{proof}
\begin{remark} \label{rere} Set $a = g'(0) >1$ and consider $z_0 = z_0(h)$ as function
of $h$. Then, after determining $z_0'(h)$ from the second equation
of (\ref{ep1}), we can use the first equation of (\ref{ep1}) to
get the following differential equation for $\epsilon_0 =
\epsilon_0(h)$:
\begin{equation}\label{de}
\epsilon_0'(h) = F(h,\epsilon_0(h)),
\end{equation}
where
$$
F(h,\epsilon) = \frac{\epsilon(2+z_0(h))}{1+h+hz_0(h)} =
\frac{\epsilon(2h\epsilon - \epsilon +0.5(h+
\sqrt{4\epsilon^2+4h^2\epsilon+h^2}))}{h(h\epsilon+0.5(h+
\sqrt{4\epsilon^2+4h^2\epsilon+h^2}))} >0.
$$
Since $F(h,\epsilon)= \epsilon/h - (1-h)/(h(1+z_0(h)h+h))$, we
find that $F(h,\epsilon) \geq \epsilon/h$ for $h \geq 1$ and
$F(h,\epsilon) \leq \epsilon/h$ for $h \in (0,1]$. Taking into
account the initial value $\epsilon_0(1) = (\ln a)^{-1}$, and
using standard comparison theory for differential equations and
the relation $\epsilon_0(0+) = 0.25/(a-1)$, we obtain that, for
every $a>1$ and $h > 0$,
\begin{equation}\label{cru} c_*(h,a) =
1/\sqrt{\epsilon_0(h,a)} \leq \min \left\{ 2\sqrt{a-1},
\sqrt{\frac{\ln a}{h}}\right\}.
\end{equation}Furthermore, $F(h,\epsilon) > 0$ so
that $c_*(h,a)$ is strictly decreasing in $h$, and we see from
(\ref{cru}) that $c_*(+\infty,a) =0$ for every $a > 1$. Now,
solving numerically the initial value problem $\epsilon_0(1) =
(\ln a)^{-1}$ for scalar differential equation (\ref{de}), we can
easily plot graphs of $c_*$ against $h$, for every fixed $g'(0)
>1$.
\end{remark}

\section{An application of Ma-Wu-Zou reduction}
\noindent Throughout this section, $\chi_{\R_-}(t)$ stands for the
indicator (or characteristic function) of $\R_-$. Next, following
the notations of Introduction and Lemma \ref{L23}, for given  $
\epsilon \in (0, \epsilon_0)$ we will denote by $\lambda_1 =
\lambda_1(\epsilon)$ and $\lambda_2 = \lambda_2(\epsilon)$ the
positive roots of $\psi(z,\epsilon) =0,$ while the roots $\epsilon
z^2 -z -1 =0$ will be denoted by $\lambda < 0 < \mu$. Notice that
$$
\lambda < 0 < \lambda_1 < \lambda_2 < \mu.
$$
Also, everywhere in this section, we require that
\begin{description}
\item[{\rm \bf(L)}] $g: \R_+ \to \R_+$ is bounded and locally
linear in some right $\delta$- neighborhood of the origin: $g(x) =
px, \ x \in [0,\delta)$, with $p > 1$, and that $g(x) \leq px$ for
all $x \geq 0$.
\end{description}
Assuming all this, we will prove the existence of non-constant
positive bounded solutions to the equation
\begin{equation}\label{twe}
\epsilon x''(t) - x'(t)-x(t)+ g(x(t-h))=0, \quad t \in \R ,
\end{equation}
satisfying $x(-\infty)=0$. Being bounded, each such solution must
satisfy the integral equation $x = Ax$:
\begin{equation}\label{iie} \hspace{0mm}
x(t) = \frac{1}{\epsilon (\mu -
\lambda)}\left\{\int_{-\infty}^te^{\lambda (t-s)}g(x(s-h))ds +
\int_t^{+\infty}e^{\mu (t-s)}g(x(s-h))ds \right\}.
\end{equation}
As it was shown by Ma, Wu and Zou \cite{hz,hz1,ma,ma1,swz,sz,wz,zw},
solving (\ref{iie}) can be successfully reduced to the determination
of fixed points of the integral operator $A$ considered in some
closed, bounded, convex and $A$- invariant subset $K$ of an
appropriate Banach space $(X, \|\cdot\|)$. In this paper, the choice
of $K \subset X$ is restricted by the following natural conditions:
(a) non-zero constant functions can not be elements of $X$; (b) the
functions $\phi^+(t) =\delta \exp (\lambda_1 t)$ and $\phi^-(t) =
\delta(e^{\lambda_1t} - e^{\lambda_2 t})\chi_{\R_-}(t) \in K$ belong
to $X$; (c) the convergence $x_n \to x$ on $K$ is equivalent to the
uniform convergence $x_n \Rightarrow x_0$ on compact subsets of
$\R$. Thus, with this in mind, for some $\rho \in (\lambda_1,\mu)$,
we  set
\begin{eqnarray*}
  X &=& \{x \in C(\R,\R): \|x\| =
\sup_{s \leq 0} e^{-\lambda_1s/2 }|x(s)|+ \sup_{s \geq 0} e^{-\rho
s}|x(s)|<
\infty\}; \\
  K &=& \{x \in X: \phi^-(t)= \delta(e^{\lambda_1 t} -
e^{\lambda_2 t})\chi_{\R_-}(t) \leq x(t) \leq \delta
e^{\lambda_1t} = \phi^+(t), \ t \in \R\}.
\end{eqnarray*}
We will consider also the following operator $L$, which can be
viewed as the 'linearization' of $A$ along the trivial steady
state:
$$
(Lx)(t) = \frac{p}{\epsilon (\mu -
\lambda)}\left\{\int_{-\infty}^te^{\lambda (t-s)}x(s-h)ds +
\int_t^{+\infty}e^{\mu (t-s)}x(s-h)ds \right\}.
$$
\begin{lemma} We have
$L\phi^{+} = \phi^{+}$. Next, for every $\psi(t) = (e^{\lambda_1
t} - e^{\nu t})\chi_{\R_-}(t) \in K$ where $\nu \in (\lambda_1,
\lambda_2]$, we have $(L\psi)(t)> \psi(t), \quad t \in \R.$
\end{lemma}
\begin{proof} It suffices to prove that $(L\psi)(t)> \psi(t)$ for $t \leq 0$. Set
$\tilde p = p\epsilon^{-1} (\mu -\lambda)^{-1}$, then we have
$$
(L\psi)(t)= \tilde p\left\{\int_{-\infty}^te^{\lambda
(t-s)}\psi(s-h)ds + \int_t^{h}e^{\mu (t-s)}\psi(s-h)ds \right\}>
$$
$$
> \tilde p\left\{\int_{-\infty}^te^{\lambda (t-s)}\psi(s-h)ds +
\int_t^{+\infty}e^{\mu (t-s)}(e^{\lambda_1(s-h)} - e^{\nu(s-h)})ds
\right\} \geq \psi(t).
$$
\end{proof}
\begin{lemma} Let assumption {\rm \bf(L)} hold and
$\epsilon \in (0, \epsilon_0)$. Then  $A(K) \subseteq K$.
\end{lemma}
\begin{proof}
We have $Ax \leq Lx \leq L\phi^+ = \phi^+$ for every $x \leq
\phi^+$. Now, if for some $s$ we have $0 <\phi^-(s) \leq x(s)$,
then $s <0$ so that $x(s) \leq \delta e^{\lambda_1s} < \delta$
implying $g(x(s)) =  px(s) \geq p\phi^-(s)$. If $\phi^-(s_1) =0$
then again $g(x(s_1)) \geq p\phi^-(s_1)=0$. Therefore $Ax \geq
L\phi^- > \phi^-$ for every $x \in K$.
\end{proof}

\begin{lemma}$K$ is a closed, bounded, convex subset of $X$
and $A:K \to K$ is completely continuos.
\end{lemma}
\begin{proof}
Notice that the convergence of a sequence in $K$ amounts to the
uniform convergence on compact subsets of $\R$. Since $g$ is a
bounded function, we have\begin{equation}\label{pro}
 |(Ax)'(t)|\leq \frac{\max_{x \geq
0}g(x)}{\epsilon (\mu - \lambda)}
\end{equation} for every $x \in K$. The
statement of this lemma follows now from the Ascoli-Arzel${\rm
\grave{a}}$ theorem combined with the Lebesgue's dominated
convergence theorem.
\end{proof}
\begin{theorem} \label{34} Assume $({\bf L})$ and let $\epsilon \in (0, \epsilon_0)$. Then the integral equation
(\ref{iie}) has a positive bounded solution in $K$.
\end{theorem}
\begin{proof}
Due to the above lemmas, we can apply the Schauder's fixed point
theorem to $A:K \to K$.
\end{proof}
\section{Bounded and uniformly persistent solutions of Eq. (\ref{twe})}
\noindent The following assumption {\rm \bf(B)} will be needed
throughout this section:
\begin{description}
\item[{\rm \bf(B)}] $g: \R_+ \to R_+$ is continuous and such that,
for some $0 <\zeta_1 < \zeta_2$,

1. $g([\zeta_1,\zeta_2])\subseteq [\zeta_1,\zeta_2]$ and
$g([0,\zeta_1])\subseteq [0,\zeta_2]$;

2. $\min_{s \in [\zeta_1,\zeta_2]}g(s) = g(\zeta_1)$;

3.  $g(x) > x$ for $x \in (0, \zeta_1]$ and $1 < g'_-(0+) \leq
g'_+(0+) < \infty$;

4. In $[0,\zeta_2]$, the equation $g(x)=x$ has exactly two fixed
points $0$ and $\kappa$.
\end{description}

\noindent In the sequel, we set $C = C([-h,0], \R)$ and $C_+ =
C([-h,0], \R_+)$. We will write
 $x(t,x'_0,\varphi)$ for the solution to the initial value problem
$x'(0) =x'_0, \ x(s) = \varphi(s)$, with $(x'_0,\varphi) \in
\R\times C$. We will also use the customary notation $x_t \in C$
where $x_t(s) = x(t+s), \ s \in [-h,0]$. In a standard way, one
parametric family of maps $S^t:\R\times C \to \R\times C, \ t \geq
0$,
$$
S^t(x'_0,\varphi) = (x'(t,x'_0,\varphi), x_t(x'_0,\varphi) )
$$
associates with autonomous equation (\ref{twe}) the continuous
semi-dynamical system $S^t$. If $x(t)$ is a bounded non negative
solution of (\ref{twe}): $0\leq x(t)\leq M_0, \ t \in \R$, then
(\ref{pro}) implies that $ |x'(t)| \leq (\epsilon (\mu -
\lambda))^{-1}\max_{x \in [0,M_0]}g(x)$. In consequence,
$(x'(t),x_t): \R \to \R \times C$ is a bounded complete trajectory
of $S^t$. Hence, with every non-negative bounded solution $x: \R \to
\R_+$ of Eq. (\ref{twe}) we can associate a compact invariant
non-empty $\omega$-limit set $\omega(x)\subset \R \times C_+$. Now,
if
$$0 \leq m = \liminf\limits_{t \to
+\infty}x(t) \leq \limsup\limits_{t \to +\infty}x(t) = M <
+\infty,$$ then, due to the compactness of $\omega(x)$, it holds
$$m=\min\limits_{\phi \in Pr\omega (x)}\min\limits_{s\in
[-h,0]}\phi(s), \quad M=\max\limits_{\phi \in Pr\omega
(x)}\max\limits_{s\in [-h,0]}\phi(s),$$ where $Pr:\R \times C_+\to
C_+$ is the canonical projection on the second factor. This
construction has the following important consequence:
\begin{lemma} \label{vg} $[m,M] \subseteq g([m,M])$.
\end{lemma}
\begin{proof}
Indeed, if there exist $M=\max_{s \in \R}x(s)$ and $m=\inf_{s \in
\R}x(s)$, a straightforward  estimation  of the right hand side of
(\ref{iie}) generates ${M} \leq \max_{{m} \leq x \leq {M}} g(x)$.
As long as the maximum $M$ is not reached, due to compactness of
$Pr\omega (x)$ we always can find a solution $z(t)$ of
$(\ref{twe})$ such that $z(0) = \max_{s \in\R}z(s) = M$ and $
\inf_{s \in\R}z(s) \geq m$. Therefore, by the above argument, ${M}
\leq \max_{{m} \leq x \leq {M}} g(x)$. The inequality ${m} \geq
\min_{{m} \leq x \leq {M}} g(x)$ can be proved in a similar way.
Thus we can conclude that  $[m,M] \subseteq g([m,M])$.
\end{proof}

\noindent From now on, we will use the following function $\xi:
\R_+ \to \R_+$,
$$
\xi(x) = \frac{\mu -\lambda}{\mu e^{-\lambda x}- \lambda e^{-\mu
x}}.
$$
\begin{lemma} \label{v} If $x$ is a solution of the boundary value problem
$x(a) = x_0, \ x'(b) =0$ for (\ref{twe}), then
$$
x(b) = \xi(b-a)\left\{x_0 + \frac{1}{\epsilon(\mu
-\lambda)}\int_{a}^b(e^{\lambda(a- u)} - e^{\mu(a-
u)})g(x(u-h))du\right\}
$$
\end{lemma}
\begin{proof} It suffices to consider the variation of constants
formula for Eq. (\ref{twe}):
$$
x(t) = Ae^{\lambda t} + Be^{\mu t} + \frac{1}{\epsilon (\mu -
\lambda)}\left\{ \int_a^te^{\lambda (t-s)}g(s)ds + \int_t^be^{\mu
(t-s)}g(s)ds \right\},
$$
where $g(s):= g(y(s-h))$.
\end{proof}
\begin{lemma} \label{41} Assume $\mathbf{(B)}$ and suppose that $\sup_{s \geq 0} g(s) \leq \zeta_2$.
Let $x\not\equiv 0$ be a non negative solution of Eq. (\ref{twe}).
Then we have the following uniform permanence property:
$$
\zeta_1  \leq \liminf\limits_{t \to +\infty} x(t) \leq
\limsup\limits_{t \to +\infty} x(t) \leq \zeta_2.
$$
\end{lemma}
\begin{proof} We follow the main lines of \cite[Corollary 12(3)]{gyt} and \cite[Theorem 3.6(a)]{LMT}.
The positivity and the upper bound of $x$ can be easily obtained
from the integral equation (\ref{iie}), so that we need only to
prove the first inequality. We proceed by contrary, supposing that
$\zeta_1
> \rho =\liminf\limits_{t \to +\infty} x(t)$. Set $\Delta =
\inf_{0\leq x\leq \zeta_1}[g(x)/x] >1$ and let $l > 0$ be large
enough and $\varepsilon > 0$ be small enough to satisfy $\Delta
(1-\varepsilon)(1-\xi(l)) > 1+\varepsilon$.   Now, we observe that
$x$ can not decrease monotonically on some semi-axis $[\tau,
+\infty)$. Indeed, in the latter case,  $x'(t) \leq 0$ and
$\zeta_1
> x(t) > x(+\infty) =0$ for sufficiently large $t$ and therefore $
x(t)\leq x(t-h) < g(x(t-h))$ (see $\mathbf{(B3)}$) implying
$$
\epsilon x''(t) = x'(t)+x(t) - g(x(t-h)) < 0.
$$
This leads us to a contradiction: $x(+\infty) = - \infty $.

In consequence,  we can indicate a local minimum point $s_0$ such
that $x(s_0) <(1+\varepsilon)\rho$ and $x(s) \geq \rho(1 -
\varepsilon)$ for all $s \in [s_0-l-h,s_0]$.  Then we have that
$g(x(s)) \geq \Delta (1-\varepsilon)\rho$ for every $s \in
[s_0-l-h,s_0]$. Indeed, if $x(s) \leq \zeta_1$, we obtain $g(x(s))
= x(s)(g(x(s))/x(s)) \geq \Delta x(s) \geq \Delta \rho(1 -
\varepsilon)$. Now, if $x(s) \in [\zeta_1,\zeta_2]$ then $g(x(s))
\geq  g(\zeta_1) \geq \Delta \zeta_1 > \Delta \rho(1 -
\varepsilon)$, see $\mathbf{(B2)}$.

Therefore, in view of Lemma \ref{v}, we have
\begin{eqnarray*}
&&\hspace{-10mm}x(s_0)= \xi(l)\left\{x(s_0-l) +
\frac{1}{\epsilon(\mu-\lambda)}\int_{s_0-l}^{s_0}(e^{\lambda(s_0-l-
u)}
- e^{\mu (s_0-l-u)})g(x(u-h))du\right\} \geq \\
   &\geq & \xi(l)\left\{
\frac{1}{\epsilon(\mu
-\lambda)}\int_{s_0-l}^{s_0}(e^{\lambda(s_0-l- u)} - e^{\mu
(s_0-l-u)})\Delta \rho(1 -
\varepsilon)du\right\} = \\
  &
= &\rho \Delta (1- \xi(l))(1 - \varepsilon) > (1+\varepsilon)\rho,
\quad {\rm a \ contradiction.}
\end{eqnarray*}
\end{proof}
\begin{lemma} \label{43} Assume $\mathbf{(B)}$ and suppose that $\sup_{s \geq 0} g(s) \leq \zeta_2$.
Let $x\not\equiv 0$ be a non negative solution of Eq. (\ref{twe})
such that $\hat{m} = \inf_{s \in \R}x(s) < \zeta_1$. Then $
\lim_{t \to -\infty} x(t) =0. $
\end{lemma}
\begin{proof} Set $\hat{M}=\sup_{s \in \R}x(s)$, Lemma \ref{vg}
guarantees that
\begin{equation}\label{si}
[\hat{m},\hat{M}] \subseteq g([\hat{m},\hat{M}]).
\end{equation}
The assumption $\mathbf{(B)}$  makes impossible the case $\hat{m}
> 0$. In consequence, $\hat{m} = 0$ and, due to Lemma \ref{41},
either $x(-\infty) = 0$ or
$$
0= \liminf\limits_{t \to -\infty} x(t) < \limsup\limits_{t \to
-\infty} x(t) = S.
$$
However, as we will show it in the continuation, the  second case
can not occur. Indeed, otherwise it is possible to indicate two
sequences of real numbers $p_n < q_n$ converging to $-\infty$ such
that $x$ has local maximum (minimum) at $p_n$ ($q_n$), $ x(q_n)<
x(s)< x(p_n)$ for all $s \in (p_n,q_n)$ and $\lim x(p_n) = S, \
\lim x(q_n) = 0$. We notice that necessarily $\lim(q_n-p_n) =
+\infty$, since in the opposite case an application of the
compactness argument leads to the following contradiction:  the
sequence of solutions $x(t-q_n)$ contains a subsequence converging
to a solution $\psi \in C(\R,\R)$ of Eq. (\ref{iie}) verifying
$\psi(0) = 0$ and $\psi(t_0) = S$, for some finite $t_0 < 0$.
Hence, $x(q_n-h) > x(q_n)$ and therefore, for sufficiently large
$n$, we get $g(x(q_n-h)) > x(q_n)$ contradicting to
$$
\epsilon x''(q_n) - x'(q_n)-x(q_n)+ g(x(q_n-h))=0,
$$
since $x''(q_n)\geq 0, \  x'(q_n) =0$.
\end{proof}
\begin{theorem} \label{mg} Assume $\mathbf{(B)}$ and suppose that
$c > c_*^+(g'_+(0+))$, then Eq. (\ref{twe}) has at least one
positive bounded solution $x(t)$ such that $$0 = \lim_{t\to
-\infty}x(t) < \zeta_1 \leq \liminf_{t\to +\infty}x(t) \leq
\sup_{s \in \R} x(s) \leq \zeta_2.$$
\end{theorem}
\begin{proof} First, let us assume that $\max_{x \geq 0} g(x) = \max_{x \in [\zeta_1,\zeta_2]}
g(x)$. Take positive integer $k$ such that $kx > g(x)$ for all $x
> 0$ and consider the following sequence
$$
g_n(u)=\left\{%
\begin{array}{ll}
    ku, & \hbox{ for} \ u \in [0,1/(nk)]; \\
    1/n, & \hbox{when }\ u \in [1/(nk), \inf g^{-1}(1/n)] \\
    g(u), & \hbox{if} \ u > \inf g^{-1}(1/n).  \\
\end{array}%
\right.
$$
of continuous functions $g_n$, all of them satisfying hypothesis
$(L)$. Obviously, $g_n$ converges uniformly to $g$ on $\R_+$. Now,
for all sufficiently large $n$, Theorem \ref{34} and Lemma
\ref{41} guarantee the existence of a positive continuous function
$x_n(t)$ such that $x_n(-\infty) =0$, $\liminf_{t \to
+\infty}x_n(t) \geq \zeta_1$, and
\begin{equation*}\label{pseea} 
x_n(t) = \frac{1}{\epsilon (\mu -
\lambda)}\left\{\int_{-\infty}^te^{\lambda (t-s)}g_n(x_n(s-h))ds +
\int_t^{+\infty}e^{\mu (t-s)}g_n(x_n(s-h))ds \right\},
\end{equation*}
Since shifted functions $x_n(s + a)$ satisfy the same integral
equation, we can assume that $x_n(0)= 0.5\zeta_1$.

Now, it is clear that the set $\{x_n\}$ is pre-compact in the
compact open topology of $C(\R,\R)$, so that we can indicate a
subsequence $x_{n_j}(t)$ which converges uniformly on compacts to
some bounded element $x \in C(\R,\R)$. Since $\lim_{j \to
+\infty}g_{n_j}(x_{n_j}(s-h)) = g(x(s-h))$ for every $s\in \R$, we
can use the Lebesgue's dominated convergence theorem to conclude
that $x$ satisfies integral equation (\ref{iie}). Finally, notice
that $x(0) = 0.5\zeta_1$ and thus $x(-\infty) = 0$ (by Lemma
\ref{43}) and $\liminf_{t \to +\infty} x(t) \geq \zeta_1$ (by
Lemma \ref{41}).

To complete the proof, we still have to analyze the case when
$\max_{x \geq 0} g(x) > \max_{x \in [\zeta_1,\zeta_2]} g(x)$.
However, this cases can be reduced to the previous one if we
redefine $g(x)$ as $g(\zeta_2)$ for all $x \geq \zeta_2$, and then
observe that $\sup_{t \in \R} x(s) \leq \zeta_2$ for every
solution obtained in the first part of the demonstration.
\end{proof}
\begin{remark} For Eq. (\ref{twe}), Theorem \ref{mg} considerably improves
the result given in Theorem 1.1 from \cite{ma1}. Indeed, the
method employed in \cite{ma1} needs essentially that
$\limsup_{u\to +0}(g'(0) - g(u)/u)u^{-\nu}$ is finite for some
$\nu \in (0,1]$ and that $g(u) < g'(0)u$.
\end{remark}
\section{Heteroclinic solutions of Eq. (\ref{twe})}
\noindent Everywhere in this section, we will assume the
hypothesis $\mathbf{(H)}$ so that all conditions of the weaker
$\mathbf{(B)}$ will be satisfied. Assume that $c > c_*^+$ and let
$x(t)$ be a bounded positive non-constant solution of Eq.
(\ref{twe}) which existence was established in Theorem \ref{mg}.
Set
$$m = \liminf\limits_{t \to
+\infty}x(t) \leq \limsup\limits_{t \to +\infty}x(t) = M.$$ Below,
we prove that $m=M=\kappa$ once $c$ is taken from the admissible
speed interval $\mathcal{C_+}$. In consequence,  the mentioned
$x(t)$ is in fact a "travelling profile" for Eq. (\ref{17}).
\begin{theorem} \label{mainex} Assume {\rm \bf(H)}. Then, for
every $c \in \mathcal{C_+}$, equation (\ref{twe}) has a positive
heteroclinic solution  $x(t)$.
\end{theorem}
\begin{proof}
 We have proved already that $\kappa, m, M \in [\zeta_1,\zeta_2]
$ and that $ [m,M]\subseteq g([m,M])$. Hence, the sequence of
intervals $\{g^k([m,M])\}$ increases and
$$[m,M]
\subseteq \cap_{k=0}^{\infty}g^k([m,M]) \subseteq
\cap_{k=0}^{\infty}g^k([\zeta_1,\zeta_2]) := [\zeta_*,\zeta^*].$$
As it can be easily checked, $g([\zeta_*,\zeta^*]) =
[\zeta_*,\zeta^*] \in \kappa$ (in fact, $[\zeta_*,\zeta^*]$ is the
global attractor of $g:[\zeta_1,\zeta_2] \to [\zeta_1,\zeta_2]$,
see \cite{gyt} for details). In particular, this means that
$g(\zeta^*) = \zeta_*$. Furthermore, there are two possibilities
for $g: [\zeta_*,\zeta^*] \to [\zeta_*,\zeta^*]$: \\ (F1) $g$ is
monotone decreasing on $[\zeta_*,\zeta^*]$, and (F2) $g$ is
unimodal on $[\zeta_*,\zeta^*]$.

Taking into account these characteristics of $g$ and applying
Lemma \ref{vg} we obtain that each of the following three
relations: $\kappa\leq x_M,$ or $\kappa \leq m\leq M,$ or $m\leq
M\leq \kappa$ implies immediately that $m=M=\kappa$. Therefore, we
will concentrate our attention only on the case when $m\leq \kappa
\leq M$ and $\kappa>x_M$ (therefore, $g'(\kappa)<0$).

Now, below we will consider only the 'unimodal' case (F2),
analogous and simpler 'monotone' case (F1) is being left to the
reader.

Let $g_1^{-1}:[g^3(x_M),\zeta^*]\to [\zeta_*,x_M]$ be the inverse
of $g$ on $[\zeta_*,x_M]$ and $\psi :[\zeta_*,\zeta^*]\to
[x_M,\zeta^*]\ni \kappa$ be the inverse of $g$ on $[x_M,
\zeta^*].$

Since $Pr\omega (x)$ is a compact invariant set, there exist a
trajectory $y(t)$ and a real $h$ such that $y(h)=M$. Since
$\kappa\leq \ M\leq \zeta^*$, we obtain from (\ref{twe}) that $
y(0)\leq \psi (M)\in [\iota, \kappa], $ where $\iota =\min
\{g_1^{-1}(M),\zeta_*\}$ if $g _1^{-1}(M)\not= \emptyset$ and
$\iota =\zeta_*$ in the opposite case.

Taking into account the boundary condition $y'(h)= 0$, setting
$g(s):= g(y(s-h))$ and then using Lemma \ref{v}, we find that
\begin{eqnarray*}
  M&=& y(h) = \xi (h)\left\{y(0) +
\frac{1}{\epsilon(\mu -\lambda)}\int_{0}^h(e^{-\lambda u}
- e^{-\mu u})g(u)du\right\} \leq \\
   &\leq & \xi (h)\left\{y(0) +
\frac{1}{\epsilon(\mu -\lambda)}\int_{0}^h(e^{-\lambda u} -
e^{-\mu u})du\max_{x \in
[m,M]}{g(x)}\right\} =\\
  &=& \xi (h) y(0) +(1- \xi (h))\max_{x \in [m,M]}{g(x)} \leq \xi \psi(M) +(1- \xi)\max_{x \in [m,M]}{g(x)}.
\end{eqnarray*}
Next, if we set $\xi = \xi(h)$ and $\theta (x)=x-\xi\psi (x)$, we
get
\begin{equation}
\theta (M)\leq (1-\xi)\max\limits_{x\in [m,M]}g(x). \label{43a}
\end{equation}

\noindent Analogously, there exists a solution $z(t)$ such that
$z(h)=m$. Then $m\in [\zeta_*,\kappa]$ and, since $z(0)\geq m >
g_1^{-1}(m),$ (when $g_1^{-1}(m) \not= \emptyset$) we have from
(\ref{twe}) that $z(0)\geq \psi (m). $ In consequence,
\begin{equation*}
m\geq \xi\psi (m)+(1-\xi)\min\limits_{x\in [m,M]}g(x),
\end{equation*}
implying
\begin{equation}
\theta (m)\geq (1-\xi)\min\limits_{x\in [m,M]} g(x). \label{46}
\end{equation}
Next, observe that
\begin{eqnarray*}
\theta ([\zeta_*,\zeta^*])=[\zeta_*-\xi\zeta^*,\zeta^*-\xi
x_M]\supseteq
[(1-\xi)\zeta_*,(1-\xi)\zeta^*]=(1-\xi)g([\zeta_*,\zeta^*])
\end{eqnarray*}
and $\theta$ is an increasing function, since $\psi$ is
decreasing. Evidently, the function $f\circ g$ is unimodal if the
functions $g, \ f$ are unimodal and increasing, respectively.
Therefore, as
$$\theta (\kappa) = \kappa - \xi\psi (\kappa) =
(1 -\xi)\kappa = (1 -\xi)g (\kappa),$$ and the function $(1
-\xi)g$ is unimodal, we have that
\[
\sigma = \theta ^{-1}\circ((1-\xi)g):[\zeta_*,\zeta^*]\rightarrow
[\zeta_*,\zeta^*]
\]
is well defined unimodal map with a unique positive fixed point
$\kappa$. Furthermore, as it was proved in \cite{gyt}, the
inequality$(Sg)(x)<0$ implies $(S\sigma)(x)<0$ and, for
$g'(\kappa)<0,$ the inequality $\left|\sigma'(\kappa)\right| \leq
1$ amounts to
\begin{equation}\label{xic}
\xi(h)
> \frac{(g'(\kappa))^{2} + g'(\kappa)}{(g'(\kappa))^{2} + 1}.
\end{equation}
From the estimates (\ref{43a}), (\ref{46}) we obtain that
$$
\lbrack m,M]\subseteq \sigma([m,M])\subseteq \sigma
^2([m,M])\subseteq ...\subseteq \sigma ^j([m,M])\subseteq ...
$$
Now, for the unimodal function $\sigma$, the last chain of
inclusions and the two conditions $|\sigma'(\kappa)| \leq 1$ and
$(S\sigma)(x) < 0$ are sufficient to deduce $m=M=\kappa$, see
Proposition \ref{SW}. Hence, Theorem \ref{mainex} and the
existence part of Theorem \ref{main} are completely proved.
\end{proof}

\begin{remark}[Nonmonotonicity]
\label{rrr} The proof of the final assertion of Theorem \ref{main}
about nonmonotonicity is in fact contained in the proof of Lemma 16
of \cite{FT} and therefore is omitted here. Now, the method applied
in \cite{FT} requires $C^2$-smoothness of $g$ at $\kappa$ and the
hyperbolicity of Eq. (\ref{cchh1}). However, as we show in a
forthcoming paper, the first condition can be weakened and second
one can be removed. In fact, it suffices to assume that $g$ is a
continuous function which is differentiable at $\kappa$. This
observation is also valid for the calculation of the minimal speed
of propagation, see Theorem \ref{ms} below.
\end{remark}
\begin{remark}[The minimal speed of propagation]
\label{rr} The minimal speed $c_*$ for reaction-diffusion
functional differential equations (\ref{17}) was already
calculated in the pioneering work \cite{sch} of K. Schaaf. See
Theorem 2.7 (i) and Lemma 2.5 of the mentioned paper (see also
\cite[Remark 4.1]{tz} for time-delayed reaction-diffusion models
different from Eq. (\ref{17})). However, the proof of minimality
of $c_*$ given in \cite{sch} seems to be incomplete. Indeed, the
proof of Proposition 2.6 in \cite{sch} uses, in an essential way,
the hyperbolicity of Eq. (\ref{char3}). See the last paragraph on
p. 591 and the formula (2.25) in \cite{sch}. Unfortunately,  for
some $c=c' < c_*$, Eq. (\ref{char3}) can be non-hyperbolic (in
other words, it can have roots on the imaginary axis) and thus
still it is necessary to disprove the existence of wavefronts
propagating at the 'critical' velocity $c'
>0$. We will do this work below, where invoking  Proposition
7.2 from \cite{FA}, we give a complete proof for the minimal speed
result.
\end{remark}
\begin{theorem} \label{ms}
If $g$ is $C^2$-smooth in some neighborhood of $x =0$,  then Eq.
(\ref{17}) does not have positive travelling wavefronts
propagating at the velocity $c < c_*$.
\end{theorem}
\begin{proof}
Suppose the theorem were false. Then there exists a travelling
profile $x: \R \to \R_+, \ x(-\infty) = 0,$ moving with positive
velocity $c' < c_*$. For the convenience of the reader, the
following part of the proof will be divided in several steps.

\noindent {\it \underline{Step (i)}} First, we suppose that $x(t)
= O(\exp (\alpha t)), \ t \to -\infty$ for some $\alpha >0$  (this
always holds when Eq. (\ref{char3}) with $c =c'$ does not have
roots $\lambda_j$ on the line $\Im \lambda=0$. See Remark
\ref{rrra} below). Next, Eq. (\ref{iie}) implies that, for $t \leq
h$,
$$
x(t) \geq \frac{1}{\epsilon (\mu - \lambda)} \int_{h}^{2h}e^{\mu
(t-s)}g(x(s-h))ds = \frac{e^{\mu t}}{\epsilon (\mu -
\lambda)}\int_{h}^{2h}e^{-\mu s}g(x(s-h))ds
$$
and therefore $x(t)$ can not decay superexponentially at
$-\infty$. In view of $C^2$-smoothness of $g$, solution $x(t)$
satisfies the linear homogeneous equation
\begin{equation}\label{eco}
\epsilon x''(t) - x'(t) - x(t) + \gamma(t)x(t-h) = 0,
\end{equation}
where $\gamma(t) =  g(x(t-h))/x(t-h) = g'(0) + O(\exp (\alpha
t))$. From \cite[Proposition 7.2]{FA}
 (see also \cite[Proposition 2.2]{hl1}), we
conclude that there are $b\geq \alpha,\ \delta>0$ and a nontrivial
eigensolution $u(t)$ of the limiting equation
\begin{equation}\label{30}
\epsilon u''(t)-u'(t)-u(t)+g'(0)u(t-h)=0
\end{equation}
on the generalized eigenspace associated with the (nonempty) set
$\Lambda$ of eigenvalues  with $\Re e\, \lambda=b$, such that
$x(t)=u(t)+O(\exp({(b+\delta)t}),\quad t\to-\infty$. On the other
hand, since $c < c_*$, we know that there are no real negative
eigenvalues of (\ref{30}): hence $\Im m\, \lambda \neq 0$ for all
$\lambda\in\Lambda$. From \cite[Lemma 2.3]{hl1}, we conclude that
$x(t)$ is oscillatory, a contradiction.

\noindent {\it \underline{Step (ii)}} Next, if $x(-\infty) = 0$
and $x(t) = O(\exp (\alpha t)), \ t \to -\infty$ does not hold for
any $\alpha>0$, then there exist a sequence $t_n \to - \infty$ and
a real number $d > 0$ such that $ x'(t_n)/x(t_n) \leq d$ for every
$n$. Indeed, otherwise for every $\alpha >0$ we can find $\tau <
0$ such that $x'(t)  \geq \alpha x(t)$ for every $t \leq \tau$.
This implies that $(x(t)\exp(-\alpha t))' \geq 0$ for $t \leq
\tau$, and therefore $x(t)\exp(-\alpha t) \leq
x(\tau)\exp(-\alpha\tau), \ t \leq \tau$. Thus $x(t) =
O(\exp(\alpha t))$ at $t = - \infty$, a contradiction. Since
$x(-\infty) = 0$, without the loss of generality we can suppose
that $x'(t_n) \geq 0$ and $x(t) \leq x(t_n)$ for all $t \leq t_n$.

Now, since $x(t)$ satisfies (\ref{eco}), we conclude that every
$y_n(t) = x(t+t_n)/x(t_n)\geq 0$ is a solution of
\begin{equation}\label{eco1}
\epsilon y''(t) - y'(t) - y(t) + \gamma(t+t_n)y(t-h) = 0,
\end{equation}
and $y_n(0)= 1 = \max_{t \leq 0}y_n(t), \ 0 \leq y_n'(0) \leq d$.
A partial integration of (\ref{eco1}) gives
\begin{equation}\label{eco2}
y'_n(t) = y'_n(0)e^{t/\epsilon} +
\frac{1}{\epsilon}\int_0^te^{(t-s)/\epsilon}(y_n(s)-\gamma(s+t_n)y_n(s-h))ds,
\end{equation}
from which we deduce the uniform boundedness of the sequence
$\{y'_n(t)\}$:
\begin{equation}\label{der2}
 |y'_n(t)| \leq d + 2g'(0), \quad t \leq 0, \ n \in \N.
\end{equation}
Together with $0 < y_n(t) \leq 1, \  t \leq 0$, inequality
(\ref{der2}) implies the compactness of  sequence $y_n(t), \ t
\leq 0$, in the compact open topology of $C(\R_-,\R_+)$.
Therefore, by the Arzel${\rm \grave{a}}$-Ascoli theorem, we can
indicate a subsequence $y_{n_j}(t)$ converging uniformly to a
continuous function $y(t)$ on bounded subsets of $\R_-$ and
additionally we may assume that $\lim_{j \to \infty} y'_{n_j}(0) =
y'_0$ exists. In consequence, $\lim_{n \to \infty}\gamma(t+t_n) =
g'(0)$ uniformly on $\R_-$. Integrating (\ref{eco2}) between $t$
and $0$ and then taking the limit as $j \to \infty$ in the
obtained expression, we establish that $y(t)$ satisfies Eq.
(\ref{30}). Additionally, $y(0)=1, y'(0) = y'_0 \in [0,d]$ and $0
\leq y(t) \leq 1, \ t \leq 0$. Moreover, we can prove  that $y(t)
> 0$ for all $ t \leq 0$.  Indeed, if $y(s) =0$ with $s < 0$ then
$y'(s)=0, y''(s) \geq 0$ so that (\ref{30}) implies $y''(s)=0$ and
$y(s-h) =0$. Setting in Lemma \ref{v} $g(x) = g'(0)x$, $a = s-h$
and $b =s$, we find that $y(t) = 0$ for all $t \in [s-2h,s-h]$.
Thus $y(t) \equiv 0$ for all $t \geq s-2h$ which is not possible
because of $y(0) =1$.

Summing up, we have proved that Eq. (\ref{30}) has a bounded
positive solution on $\R_-$.  As it was established in \cite[Lemma
A.1]{hl1},  this solution does not decay superexponentially. In
consequence, arguing as in the last paragraphs of {\it {Step
(i)}}, we finalize the proof of {\it {Step (ii)}} and Theorem
\ref{ms}.
\end{proof}

\begin{remark}[Asymptotic formulae]
\label{rrra} Under conditions $c > c_*$ as well as of
$C^2$-smoothness of $g$ and hyperbolicity of (\ref{twe}) at zero, we
can repeat {\it {Step (i)}} above to derive asymptotic formulae for
heteroclinic solutions. Indeed, $x(-\infty) =0$ together with
(\ref{iie}) implies $x'(-\infty) =0$. Hence, the solution
$(x(t),x'(t))$ belongs to the unstable manifold of the trivial
equilibrium of Eq. (\ref{twe}). Thus $|x(t)|+ |x'(t)| = O(e^{\alpha
t}), \ t \to - \infty,$ for some $\alpha >0$ (see \cite[Section
VIII.4]{dglw}). On the other hand, we know that $x(t)$ is not
decaying superexponentially at $-\infty$. Since $x(t)$ satisfies the
homogeneous equation (\ref{eco}), due to \cite[Proposition 7.2]{FA},
we get, for some $k_1$ and every small $\delta >0$, that $x(t) =
k_1e^{\lambda_1 t} + O(e^{(\lambda_1 + \delta) t}) = O(e^{\lambda_1
t}),  \ t \to -\infty$. In consequence, $x(t)$  satisfies the
inhomogeneous equation
$$
\epsilon x''(t) - x'(t)-x(t)+ g'(0)x(t-h)= H(t), \quad t \in \R,
$$
where $H(t) = g'(0)x(t-h) - g(x(t-h)) = O(\exp(2\lambda_1 t))$ as $t
\to - \infty$. Now, from \cite[Proposition 7.1]{FA} and Lemma
\ref{L23} we deduce that there exist $k_1,k_2, \ k_1^2 + k_2^2
\not=0,$ such that, for every fixed $\delta > 0$,
$$
x(t)=\left\{%
\begin{array}{lll}
    k_1e^{\lambda_1 t} + k_2e^{\lambda_2 t} + O(e^{(2\lambda_1- \delta) t}),
    \ t \to -\infty, & \hbox{if}\ 2\lambda_1 > \lambda_2, \ k_1\not=0 ; \\
    k_1e^{\lambda_1 t} +  O(e^{(2\lambda_1- \delta) t}), \ t \to -\infty, & \hbox{if}\ 2\lambda_1 \leq \lambda_2, \
    k_1\not=0;\\
     k_2e^{\lambda_2 t} +  O(e^{(2\lambda_2- \delta) t}), \ t \to -\infty, & \hbox{if} \
    k_1=0, \ k_2\not=0.
\end{array}%
\right.
$$
\noindent Hence, taking an appropriate $t_0$ and $i \in \{1,2\}$,
we obtain the following asymptotic formula for travelling
wavefront profiles:
\begin{equation}\label{afo}
x(t-t_0) = \exp (\lambda_i t) + O(\exp((2\lambda_i- \delta) t)),
\quad t \to - \infty.
\end{equation}
To get a similar formula for $x'(t-t_0)$, it suffices to
differentiate (\ref{iie}) and then use (\ref{afo}) in the right
hand side of the obtained integral representation for $x'(t-t_0)$,
see \cite[Theorem 14]{FT} for details.
\end{remark}
\noindent Finally, in Proposition \ref{co} below we prove that
$c^* \leq c^*_{opt}$. In this way we show that the existence
result of Theorem \ref{main} is consistent with the conjecture
about the exact speed interval for Eq. (\ref{17}) that was
proposed in Introduction.
\begin{proposition} \label{co} Let us assume (\ref{xic}) and set $a = g'(\kappa)$. Then the
characteristic equation (\ref{char3}) has only one (real) root
$\lambda_0$ in the right half plane $\Re z \geq 0$.
\end{proposition}
\begin{proof} On the contrary, let us suppose that Eq. (\ref{char3}) has a pair $\lambda_{1,2}$
of complex conjugate roots with $\Re \lambda_{1,2}\geq 0$ for some
fixed values $h_0,\  g'(\kappa)$ and $\epsilon_0$. Notice that, by
virtue of Lemma \ref{c1}, in this case
$$\xi = \xi(h_0) >
(g'(\kappa))^{2} + g'(\kappa))/((g'(\kappa))^{2} + 1) \geq
e^{-h_0}.
$$
Here, we use the fact that (\ref{gsc}) implies $\Re \lambda_j < 0$
for every solution $\lambda_j$ of the characteristic equation
$\phi(\lambda) =0$. See \cite{gyt,ltt}.

Now, keeping $\epsilon_0$ and $g'(\kappa)$ fixed and considering
$h \in (0,h_0]$, we define $h_1 \leq h_2$ as the minimal value of
$h$ for which equation (\ref{char3}) has roots $\mu = \mu(h_1) =
\pm i \omega$ on the imaginary axis. Since
$$
\Re \frac{d\mu(h_1)}{dh} = \omega^2(1 +2\epsilon +2
\epsilon^2\omega^2)/((1+h(\epsilon \omega^2
+1))^2+\omega^2(2\epsilon-h)^2) >0,
$$
by the Hopf bifurcation theorem \cite{hale}, Eq. (\ref{twe}) has a
supercritical Hopf bifurcation at $h=h_1$. Therefore, for all $h
\in (h_1,h_0)$ close to $h_1$,  equation (\ref{twe}) has a stable
periodic solution bifurcating from the equilibrium $\kappa$. By
Theorem \ref{main}, this implies immediately that
$$
\xi(h_1) \leq \frac{(g'(\kappa))^{2} +
g'(\kappa)}{(g'(\kappa))^{2} + 1} < \xi(h_0),
$$
a contradiction since $\xi =\xi(h)$ is a strictly decreasing
function.
\end{proof}
\section*{Acknowledgments}
\quad  This work was partially supported by  FONDECYT (Chile) and by
the University of Talca, program ``Reticulados y Ecuaciones" (Sergei
Trofimchuk).

\end{document}